\documentclass[11pt]{article}
\usepackage{epsfig}
\usepackage{amssymb}
\usepackage{amsmath}
\usepackage{mathrsfs}
\newcommand{\be}{\begin{equation}}
\newcommand{\ee}{\end{equation}}
\newcommand{\bea}{\begin{eqnarray}}
\newcommand{\eea}{\end{eqnarray}}

\usepackage{bm}
\usepackage{color}
\usepackage{graphicx}
\usepackage{cite}

\newcommand{\red}[1]{{\color{red}#1}}

\begin{document}


\title{A New Class of Linear Relations for Scalar Partitions}
\author{Boris Y. Rubinstein\\
Stowers Institute for Medical Research
\\1000 50th St., Kansas City, MO 64110, U.S.A.}
\date{\today}

\maketitle
\begin{abstract}
A scalar integer partition problem asks for a number of nonnegative
integer solutions to a linear Diophantine
equation with integer positive coefficients. 
The manuscript discusses an algorithm of derivation of 
novel linear relations involving the finite number of scalar partitions.
The algorithm employs the Cayley theorem about 
the reduction of a double partition to a sum of scalar partitions
based on the variable elimination procedure.
\end{abstract}





\section{Integer scalar partitions}
\label{intro0}

The problem of integer partition into a set of integers 
is equivalent to counting number
of nonnegative integer solutions of the Diophantine equation
\be
s = \sum_{i=1}^m d_i x_i = {\bf d}\cdot{\bf x}.
\label{coin1}
\ee
A scalar partition function $W(s,{\bf d}) \equiv W(s,\{d_1,d_2,\ldots,d_m\})$
solving the above problem is a
number of partitions of an integer $s$ into positive integer generators
$\{d_1,d_2,\ldots,d_m\}$. 
The generating function for $W(s,{\bf d})$
has a form
\be
G(t,{\bf d})=\prod_{i=1}^m\frac{1}{1-t^{d_{i}}}
 =\sum_{s=0}^{\infty} W(s,{\bf d})\;t^s\;.
\label{WGF}
\ee
Cayley discovered \cite{Cayley1855} 
the splitting of the 
scalar partition into periodic and non-periodic parts and and later Sylvester 
showed that it
might be presented as a sum of "waves"
\be
W(s,{\bf d}) = \sum_{j=1} W_j(s,{\bf d})\;,
\label{SylvWavesExpand}
\ee
where the summation runs over all distinct factors
of the elements $d_i$ of the generator vector ${\bf d}$, and each
wave $W_j(s,{\bf d})$ is a quasipolynomial in $s$
closely related to prime roots $x=\rho_j$ 
of the equation $1-x^j=0$.
Each Sylvester wave $W_j(s,{\bf d})$ and thus the scalar partition $W(s,{\bf d})$ 
can be expressed 
through 
the Bernoulli polynomials of higher order \cite{Rub04}.

The definition (\ref{WGF}) of the generating function for $W(s,{\bf d})$ implies the recursion relation
\be
W(s,{\bf d}) - W(s-d_i,{\bf d}) = W(s,\{d_{1},d_{2},\ldots,d_{i-1},d_{i+1},\ldots,d_{m}\}),
\label{SPFrecursion}
\ee
involving three scalar partitions -- two with $m$ generators and
one with $m-1$ generators. 
As the selection of the 
generator $d_i$ to be dropped from the vector ${\bf d}$ is arbitrary, one can produce 
$m$ relations similar to (\ref{SPFrecursion}). 

There are just a few relations involving the partitions $W(s,{\bf d})$ 
are known in addition to (\ref{SPFrecursion}).
For example it was shown in \cite{RubFel03} that the similar relation holds for each Sylvester wave $W_j(s,{\bf d})$ 
\be
W_j(s,{\bf d}) - W_j(s-d_i,{\bf d}) = W_j(s,\{d_{1},d_{2},\ldots,d_{i-1},d_{i+1},\ldots,d_{m}\}),
\label{SPFrecursion_wave}
\ee
In his pioneering work \cite{Cayley1855} on the analytical expression for the scalar partitions Cayley
noted that the function $W(s,{\bf d})$ satisfies the 
following parity property for the positive values of $s$
\be
W(s,{\bf d}) = (-1)^{m+1} W(-s-\sigma_1({\bf d}),{\bf d}),
\quad \quad
\sigma_1({\bf d}) = \sum_{i=1}^m d_i,
\label{SPF_symm}
\ee
and pointed out that $W(s,{\bf d})$ vanishes at all integer negative points in the range
$-\sigma_1({\bf d}) < s < 0$. It is worth to note that Cayley 
considered the parity property as ``\ldots uninterpretable  in  the  theory  of  partitions'' \cite{Cayley1855}.
The same relation was derived independently \cite{FelRub01} using 
only the resursion (\ref{SPFrecursion}).


The manuscript introduces the linear relations of different structure -- (a)
each term in the relation has the same number $m$ of generators; (b)
the number of terms in the relation is $(m+1)$ -- one more than the number of generators;
and (c) the number of such relations is unlimited.
The algorithm of derivation of these relations is based on the reduction of 
a specially constructed double partition into a sum of scalar partitions.
The general idea of the reduction was proposed by Sylvester in  \cite{Sylv1}
and the actual reduction algorithm is due to Cayley \cite{Cayley1860} who also
discussed the limitations of the method.
The manuscript shows how these conditions should be employed for 
the double partition construction to obtain the required relations.


\section{Sylvester-Cayley reduction algorithm for double partitions}
\label{SylvCayley}
Double partitions used for the derivation of the scalar partition relations
represent the simplest 
case of {\it vector partitions}.  
The vector partition function $W({\bf s},{\bf D})$ counts the number of integer
nonnegative solutions ${\bf x} \ge 0$
to a linear system ${\bf s} = {\bf D} \cdot {\bf x}$, where
${\bf D}= \{{\bf c}_1,{\bf c}_2,\ldots,{\bf c}_m\}$ is a nonnegative integer $l \times m$ generator matrix ($l < m$) 
made of columns ${\bf c}_i =\{c_{i1},c_{i2},\ldots,c_{il}\}^T,\ (1 \le i \le m)$ 
with $\{\cdot\}^T$ denoting transposition of a vector.

Sylvester suggested an iterative procedure of 
reduction of a vector partition into a sum of scalar partitions \cite{Sylv1}; the procedure is based on the 
variable elimination from the system of linear equations ${\bf s} = {\bf D} \cdot {\bf x}$.
Based on this approach Cayley 
developed  an algorithm of double partition reduction
and established a set of conditions of the method applicability \cite{Cayley1860}.
Cayley showed that each column ${\bf c}_i$ of the two-row matrix ${\bf D}$ gives rise 
to a scalar partition and elements of this column must be relatively prime (we call it P-column). 
Another limitation of the Sylvester-Cayley method is that the columns ${\bf c}_i$
should represent {\it noncollinear} vectors. 
The practical application of the Cayley algorithm to the computation of the  Gaussian polynomial
coefficients is considered in \cite{RubSylvCayl}. 
The details of the Cayley reduction method for double partitions and its 
modification that works without any limitations for arbitrary double partition 
are described in \cite{RubDouble2023}. 

It should be noted that 
similar restictions to those discussed by Cayley in \cite{Cayley1860} were earlier 
mentioned in \cite{Sylv1} for vector partitions with arbitrary dimension $l$ of vector ${\bf s}$.
Later Sylvester in his lectures \cite{SylvLectures1859} offered an alternative approach to lifting
these limitations and it appears to be quite useful for our derivation procedure.

Consider the Sylvester idea in more details. The variable elimination algorithm for 
vector partition reduction fails when the generator matrix ${\bf D}$ has
either collinear columns or nonprime columns (NP-columns) or both, where
the elements of the NP-column ${\bf c}$ have the greatest common divisor (GCD)
of its elements $\mbox{gcd}({\bf c})>1$.
Sylvester suggested to add  to the system ${\bf s} = {\bf D} \cdot {\bf x}$ a single auxiliary equation
and a new unknown $x_{m+1}$ 
in order to convert the original system into a larger one 
$\hat{\bf s} = \hat {\bf D} \cdot \hat {\bf x}$,
where $\hat {\bf s} = \{s_0,s_1,s_2,\ldots,s_l\}$ and $\hat {\bf x} = \{x_1,x_2,\ldots,x_m,x_{m+1}\}$.
The matrix 
$\hat {\bf D}=\{\hat {\bf c}_1,\hat {\bf c}_2,\ldots,\hat {\bf c}_m,\hat {\bf c}_{m+1}\}$ 
consists of $(m+1)$ columns $\hat {\bf c}_{i} =\{c_{i0},c_{i1},c_{i2},\ldots,c_{il}\}^T$.
The elements $c_{i0}$ and the value of $s_0$ in the added equation must be chosen to satisfy two conditions -- 
(a) each nonnegative integer solution of the original system ${\bf s} ={\bf D} \cdot {\bf x}$
should correspond to a single solution of the system $\hat {\bf s} = \hat {\bf D} \cdot \hat {\bf x}$,
and (b) the matrix  $\hat {\bf D}$ should be void of the collinear and NP-columns.
The first condition leads to 
\be
W(\hat{\bf s},\hat{\bf D}) = W({\bf s},{\bf D}),
\label{solution_number}
\ee
while the second one guarantees that the vector partition $W(\hat{\bf s},\hat{\bf D})$ admits the 
variable elimination 
procedure of reduction (or at least the first step of it).


\section{Linear relations for scalar partitions}
\label{LinRel}
Consider an application of the Sylvester method to the scalar partition $W(s,{\bf d})$ 
and add to the single original equation $s = {\bf d} \cdot {\bf x}$ an equation 
$\sigma = x_{m+1} - \bm \delta \cdot {\bf x}$, producing the system of two equations
\be
\sigma =  x_{m+1} - \sum_{i=1}^{m} \delta_i x_i,
\quad
s = \sum_{i=1}^{m} d_i x_i,
\label{system2}
\ee
with positive $\delta_i$ and nonnegative $\sigma$ (see the fifth lecture in\cite{SylvLectures1859}).
This system
corresponds to the double partition $W(\hat{\bf s},\hat{\bf D})$ with 
\be
\hat{\bf D} =
\left[
\begin{array}{ccccc}
-\delta_1 & -\delta_2 & \ldots &-\delta_m &1   \\
d_1 & d_2 & \ldots & d_m & 0 
\end{array}
\right],
\quad
\hat{\bf s} =
\left[
\begin{array}{c}
\sigma  \\
s
\end{array}
\right].
\label{matr_hatD}
\ee
Each integer nonnegative solution ${\bf x}=\bm\xi =\{\xi_1,\xi_2,\ldots,\xi_m\}$
of the second equation $s = {\bf d} \cdot {\bf x}$
corresponds to the following solution of (\ref{system2}) 
\be
x_i = \xi_i, 
\quad
(1 \le i \le m),
\quad
x_{m+1} = \sigma +  \sum_{i=1}^{m} \delta_i \xi_i.
\label{system2sol}
\ee
We observe that the condition (\ref{solution_number}) is satisfied as for each solution
$\bm\xi$ of the original Diophantine equation $s = {\bf d} \cdot {\bf x}$ there 
exists the single solution (\ref{system2sol}) of the system (\ref{system2}).

The process of the double partition reduction through the 
elimination of the variables $x_i$ one by one is discussed 
in \cite{Cayley1860,RubSylvCayl,RubDouble2023} for the nonnegative matrix ${\bf D}$ and 
its generalization to the arbitrary matrices is given in \cite{Rub2025}. 
As each column of $\hat{\bf D}$ leads to a single scalar partition the result of the reduction is a sum of $(m+1)$ scalar partitions.
It is shown however in \cite{Rub2025} that the contribution of the last column $\hat {\bf c}_{m+1}=\{1,0\}^T$
is zero and 
we end up with the scalar partition $W(s,{\bf d})$ expressed as a sum of 
$m$ scalar partitions $W(s_i,{\bf d}_i)$ 
with ${\bf d}_i=\{d_{i,1},d_{i,2},\ldots,d_{i,m},d_{i,m+1}\}$, that satisfy the
linear relation 
\be
w(s,{\bf d},\bm \delta) \equiv W(s,{\bf d}) - \sum_{i=1}^m W(s_i,{\bf d}_i) = 0,
\label{lin_rel} \\
\ee
where
\be
d_{i,i} = d_i,
\quad
d_{i,j} = \delta_i d_j - \delta_j d_i, \ (1\le j\le m, j\ne i),
\quad
s_i = s \delta_i + \sigma d_i.
\label{part_scalar}
\ee
As for $i \ne j$ we have $d_{i,j}+d_{j,i} = 0$, the quantities $\sigma_1({\bf d}_i)$ 
sum up to $\sigma_1({\bf d})$ and we obtain the numerical relation 
\be
\sum_{i=1}^m \sigma_1({\bf d}_i) = \sigma_1({\bf d})
\quad
\Rightarrow
\quad
\frac{1}{\sigma_1({\bf d})} \sum_{i=1}^m \sigma_1({\bf d}_i) = 1.
\label{sigma_di_rel}
\ee

Consider the term $W(s_i,{\bf d}_i)$ -- in the number theory context it is defined only for nonnegative integer values
of $s_i$. As $s_i \ge s \delta_i \ge 0$ for positive $\delta_i$ and nonnegative $s$, without loss of generality
we can set $\sigma=0$ and further we always use  $s_i = s \delta_i$.
As all columns $\hat {\bf c}_{j}=\{-\delta_j,d_j\}^T$ 
are noncollinear, 
all elements of the vector ${\bf d}_i$ are nonzero but some of them might be negative
(say, $d_{i,j_k}<0$ for $1 \le k \le K_i$).
Noting that 
$(1-t^{-a})^{-1} = -t^{a}(1-t^{a})^{-1}$ we find an equivalent scalar partition with positive generators only
\be
W(s_i,{\bf d}_i) =
(-1)^{K_i} W(s \delta_i + \sum_{k=1}^{K_i} d_{i,j_k},|{\bf d}_i|),
\quad
|{\bf d}_i| = \{|d_{ij}|\}.
\label{scalar_pos_gen}
\ee

The selection of the positive vector $\bm \delta$ is a (partially heuristic) process and we discuss it
in two major cases -- (a) all $d_i$ in the set ${\bf d}$ are unique,
and (b) the set ${\bf d}$ contains duplicate generators.

In the case of the unique generators 
the simplest choice for $\bm \delta$ is to set {\it all} its elements $\delta_i=1$ to unity, i.e., $\bm \delta = {\bf 1},$ 
and we obtain 
\be
\hat{\bf D} =
\left [
\begin{array}{ccccc}
-1 & -1 & \ldots &-1 &1   \\
d_1 & d_2 & \ldots & d_m & 0 
\end{array}
\right ].
\label{matr_hatD_unique}
\ee
From (\ref{part_scalar}) we find the general term $W(s_i,{\bf d}_i)$ with
\be
s_i=s,
\quad d_{i,i} = d_i,
\quad
d_{i,j} = d_j-d_i, \ (1\le j\le m, j\ne i).
\label{part_scalar1}
\ee
We observe that this choice of $\delta_i$ makes all 
terms in the relation (\ref{lin_rel}) to have the same argument $s$.
The example of the relation with $\bm \delta = {\bf 1}$ is discussed in Section \ref{Unique}.

Other possibilities assume one or more $\delta_i > 1$ and these values should be taken
with caution not to produce a NP-column or/and a group of collinear columns.
For example, $\delta_i$ should be odd for any even $d_i$, and 
if $k \mid d_{r}$ one should set $\delta_r \ne k$.


When the original scalar partition has several duplicate elements in the set ${\bf d}$, 
say $d_1=d_2=\ldots =d_k$ one has to choose all different $\delta_i,\ 1 \le i \le k,$ such that also
$\mbox{gcd}(d_i,\delta_i) = 1$; we consider the representative example in Section \ref{Nonunique}.


\section{Linear relations analysis}
\label{properties}
The relation (\ref{SylvWavesExpand}) presents
the scalar  partition $W(s,{\bf d})$ as a sum of the pure polynomial 
$W_1(s,{\bf d})$ and several quasipolynomials $W_j(s,{\bf d})$, each being
a superposition of polynomials multiplied by the periodic function of period $j$ \cite{Sylv2,Rub04}. 
The order $n_1=m-1$ of the polynomial part $W_1(s,{\bf d})$ is one less the number $m$ of 
generators in ${\bf d}$, while the order $n_j=m_j-1$ of the polynomial factor in $W_j(s,{\bf d})$
is one less the multiplicity $m_j$ of the factor $j$ among all generators $d_i$.
The function $w(s,{\bf d},\bm \delta)$ in (\ref{lin_rel}) has $m+1$  scalar partition terms each one having the 
polynomial part $W_1(s \delta_i,{\bf d}_i)$ of order $m-1$ as well as multiple 
quasipolynomials with polynomial factors of the order $0 \le m_j < m-1$ of the periodic functions with integer period $j>1$. 

\subsection{Analytical properties}

It should be underlined that in the integer partition context the relation (\ref{lin_rel}) 
is satisfied for all nonnegative {\it integer} values of $s$.
As the individual summand in (\ref{lin_rel}) is represented by a quasipolynomial that has periodic
factors with integer period one may extend $W(s,{\bf d})$ to a function of continuous argument \cite{Rub04} by
choosing the $j$-periodic function $\psi_j(s)$ to be 
the prime radical circulator introduced by Cayley in \cite{Cayley1855}
\be
\psi_j(s) = \Psi_j(s) = \sum_{\rho_j} \rho_j^s, 
\quad
\rho_j = \exp(2\pi i n/j),
\label{Psi_def}
\ee
where the summation is made over all primitive roots of unity $\rho_j$ with $n$
relatively prime to $j$ (including unity) and smaller than $j$. 


The function $w(s,{\bf d},\bm \delta)$ has a quite unusual behavior -- it vanishes in {\it all integer points} of the real axis
while its behavior between them requires further analysis.
The general expression for the factor $U_k(s,\bm \delta)$ of $s^k, \ (0 \le k \le m-1)$ 
in the quasipolynomial $w(s,{\bf d},\bm \delta)$ 
has a constant term $A_k$ and the finite number of $\Psi_j(s  \delta_i - r_j)$ terms oscillating around zero:
\be
U_k(s,\bm \delta) = A_k +\sum_{i=0}^{m}  \sum_{j>1}\ \sum_{r_j=1}^{j} C_{k,j,r_j} \Psi_j(s \delta_i - r_j),
 \quad 0 \le k \le m-1,
\label{w_general_term}
\ee
where we introduce the notation $\delta_0 \equiv 1$.
The value of $U_k(s,\bm \delta)$ should be zero at the infinite number of integer values of $s$ and it is possible only 
when $A_k = 0$. 
Then the purely polynomial part $w_1(s,{\bf d},\bm \delta)$ 
defined via 
the corresponding part of $w(s,{\bf d},\bm \delta)$ summands must vanish at any real $s$
\be
w_1(s,{\bf d},\bm \delta) \equiv W_1(s,{\bf d}) - \sum_{i=1}^m W_1(s \delta_i,{\bf d}_i) = \sum_{k=0}^{m-1} A_k s^k 
= 0.
\label{poly_lin_rel0}
\ee
Returning to the 
expression of $w(s,{\bf d},\bm \delta)$  
write it as a multiple sum of powers of $s$ multiplied by the periodic triginometric functions
$$
w(s,{\bf d},\bm \delta) = \sum_{i=0}^{m} \sum_{k=0}^{m-2} 
\sum_{j>1}\ \sum_{r_j=1}^{j} C_{k,j,r_j}  s^k  \Psi_j(s \delta_i - r_j).
$$
It should be noted that while $\Psi_j(s - r_j)$ has integer period $j$,  the period of the function
$\Psi_j(s \delta_i - r_j)$ for $\delta_i>1$ is given in the general case by a rational fraction $j/\delta_i$.
This observation makes $\bm \delta = {\bf 1}$ to be the special case as 
then all terms in $w(s,{\bf d},{\bf 1})$ have the integer periods only. The  number of these terms is finite
and as they cancel each other at infinite number of integer points they should cancel identically for
both integer and real argument values.

On the contrary, when $\bm \delta \ne {\bf 1}$ the complete term cancelling is possible at the integer points only
but the terms with the rational periods $j/\delta_i$ in $W(s \delta_i,{\bf d}_i)$ 
cannot compensate the integer periods terms in $W(s,{\bf d})$ for the real $s$ values.
As the result the function $w(s,{\bf d},\bm \delta)$ between the integers points is not identically zero.

\subsection{Polynomial part}

The polynomial part $W_1(s,{\bf d})$ in (\ref{poly_lin_rel0}) has the following functional representation
\be
W_1(s,{\bf d}) = \frac{1}{(m-1)! \pi({\bf d})} B^{(m)}_{m-1}(s+\sigma_1({\bf d}),{\bf d}), 
\quad
\quad
\pi({\bf d}) = \prod_{i=1}^m d_i,
\label{poly_part_def}
\ee
and the Bernoulli polynomials of higher order $B^{(m)}_{n}(s,{\bf d})$ are defined by the 
generating function \cite{Norlund1924} 
$$
\frac{e^{st} \prod_{i=1}^m d_i}{\prod_{i=1}^m (e^{d_it}-1)} = 
\sum_{n=0}^{\infty} B^{(m)}_{n}(s,{\bf d})\frac{t^n}{n!}.
$$
Substitution of (\ref{poly_part_def}) into (\ref{poly_lin_rel0}) produces the linear relation for
the Bernoulli polynomials of higher order
\be
\frac{B^{(m)}_{m-1}(s+\sigma_1({\bf d}),{\bf d})}{\pi({\bf d})}  - 
\sum_{i=1}^m  \frac{B^{(m)}_{m-1}(s \delta_i+\sigma_1({\bf d}_i),{\bf d}_i)}{\pi({\bf d}_i)}  = 0.
\label{poly_lin_rel}
\ee
Considering the above formula first note that 
$
B^{(m)}_{n}(s+\sigma_1({\bf d}),{\bf d}) = B^{(m)}_{n}(s,-{\bf d})= (-1)^n B^{(m)}_{n}({\bf d}),
$
and use the binomial relation \cite{Norlund1924}
to write it as
\be
B^{(m)}_{n}(s,{\bf d}) = 
\sum_{k=0}^{n} \binom{n}{k} s^k B^{(m)}_{n-k}({\bf d}),
\label{bernoulli_binom}
\ee
where $B^{(m)}_{n}({\bf d})$ denotes the Bernoulli number of higher order.
The relations (\ref{poly_lin_rel}) and (\ref{bernoulli_binom}) lead to
\be
\frac{B^{(m)}_{m-k-1}({\bf d})}{\pi({\bf d})} - \sum_{i=1}^m   
\frac{\delta_i^k B^{(m)}_{m-k-1}({\bf d}_i)}{\pi({\bf d}_i)} = 0.
\label{Bernoulli_lin_rel}
\ee

\subsection{Relations for complete Bell polynomials}

In \cite{Rub09} the following relation was established between the Bernoulli polynomials of higher order
and the complete Bell polynomials
\be
B^{(m)}_n(s,{\bf d}) = \mathbb B_n(s+a_1,a_2,\ldots),
\quad
a_r = (-1)^{r-1}B_r \sigma_r({\bf d})/r,
\quad
\sigma_r({\bf d})=\sum_{i=1}^{m}d_i^r,
\label{BernoulliBell}
\ee
where $\sigma_r({\bf d})$ denotes a power sum of the generators ${\bf d}$ and
the complete Bell polynomials $ \mathbb B_n(a_1,a_2,\ldots)$ are defined by the generating function \cite{Riordan}
\be
\exp \left(\sum_{i=1}^{\infty}\frac{a_i}{i!}t^i \right) =
\sum_{i=0}^{\infty} \frac{\mathbb B_i(a_1,a_2,\ldots)}{i!}t^i.
\label{BellGF}
\ee
The expression for $ \mathbb B_{k}({\bf a})\equiv \mathbb B_{k}(a_1,a_2,\ldots)$ depends
on the first $k$ elements of the vector ${\bf a}$ only.
Write down the explicit expressions for $ \mathbb B_{k}({\bf a})$ for small $k\le 5$:
\bea
&& \mathbb B_{0}({\bf a}) = 1,
\quad\quad
 \mathbb B_{1}({\bf a}) = a_1,
\quad\quad
\mathbb B_{2}({\bf a}) =  a_1^2+ a_2,
\quad
\mathbb B_{3}({\bf a}) =  a_1^3+ 3a_1a_2, 
\nonumber\\
&&  \mathbb B_{4}({\bf a}) =  a_1^4+ 6a_1^2a_2+ 3a_2^2+a_4, 
\quad
\mathbb B_{5}({\bf a}) =  a_1^5+ 10a_1^3a_2 + 15a_1a_2^2 
+ 5a_1a_4 
,
 \label{Bell_explicit} 
\eea
where we take into account that $B_{2k+1}=a_{2k+1}=0$ for $k \ge 1$. 
Use (\ref{BernoulliBell}) in the relation (\ref{Bernoulli_lin_rel}) to obtain
\be
\frac{\mathbb B_{m-k-1}(a_1,a_2,\ldots)}{\pi({\bf d})} - \sum_{i=1}^m   
\frac{\delta_i^k  \mathbb B_{m-k-1}(a_{i,1},a_{i,2},\ldots)}{\pi({\bf d}_i)} = 0,
\quad
a_{i,r} =(-1)^{r-1}\sigma_r({\bf d}_i)/r.
\label{Bell_lin_rel}
\ee
This result together with (\ref{Bell_explicit}) allows to derive a number of 
numerical relations in addition to (\ref{sigma_di_rel}).
Use $k=m-1$ and $k=m-2$ in (\ref{Bell_lin_rel}) and find
\be 
\pi({\bf d})  \;\sum_{i=1}^m  \frac{\delta_i^{m-1}}{\pi({\bf d}_i)} = 1,
\quad
\frac{\pi({\bf d})}{\sigma_1({\bf d})}  \;\sum_{i=1}^m  \frac{\delta_i^{m-2}\sigma_1({\bf d}_i)}{\pi({\bf d}_i)} = 1.
\label{poly_lead_terms_rel}
\ee
With $k=m-3$ we find $ \mathbb B_{2}({\bf a}) =  a_1^2+ a_2 = \sigma_1^2({\bf d})/4 -\sigma_2({\bf d})/12$
and obtain
\be 
\frac{\pi({\bf d})}{3\sigma_1^2({\bf d}) -\sigma_2({\bf d})}  
\;\sum_{i=1}^m  \delta_i^{m-3}\;\frac{3\sigma_1^2({\bf d}_i) -\sigma_2({\bf d}_i)}{\pi({\bf d}_i)} = 1.
\label{poly_3rd_term_rel}
\ee
For each set $\bm \delta$ the total number of the numerical relations (\ref{sigma_di_rel}) and (\ref{Bell_lin_rel}) involving the products and 
power sums of the integer generators ${\bf d}$ and ${\bf d}_i$ is equal to $m$. 

The relation (\ref{Bell_lin_rel}) leads us to a conjecture 
for arbitrary 
sets ${\bf x} = \{x_1,x_2,\ldots,x_m\}, {\bf y} = \{y_1,y_2,\ldots,y_m\}$ of $m$ variables each 
\be
\frac{B^{(m)}_{m-k-1}({\bf x})}{\pi({\bf x})} - \sum_{i=1}^m   
\frac{y_i^{k} B^{(m)}_{m-k-1}({\bf s}_i)}{\pi({\bf s}_i)} = 0, 
\quad k<m,
\quad 
s_{i,j} = y_i x_j - y_j x_i + x_i \delta_{ij}, \ (1\le j\le m),
\label{Bell_lin_polynomials} 
\ee
where $\delta_{ij}$ denotes the Kronecker delta symbol.

\section{Numerical examples}
\label{examples}

In this section  
we discusss several examples illustrating the general results established above.

\subsection{Unique generators}
\label{Unique}

\begin{figure}[h!]
\begin{center}
\begin{tabular}{cc}
\includegraphics[height=4.8cm]{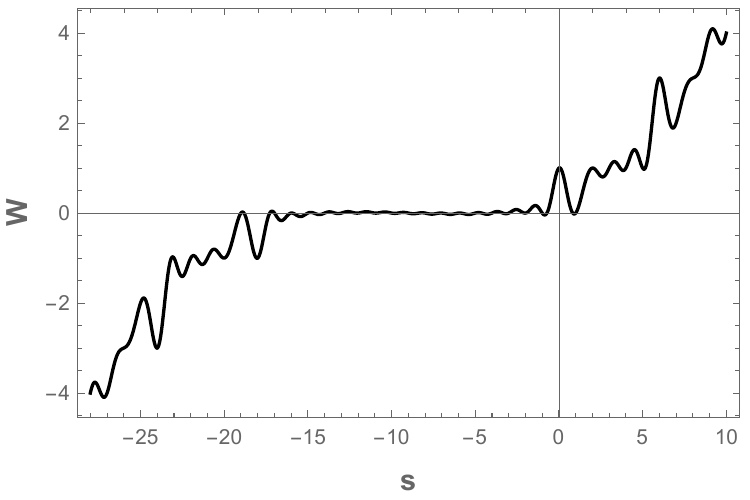}&
\includegraphics[height=4.8cm]{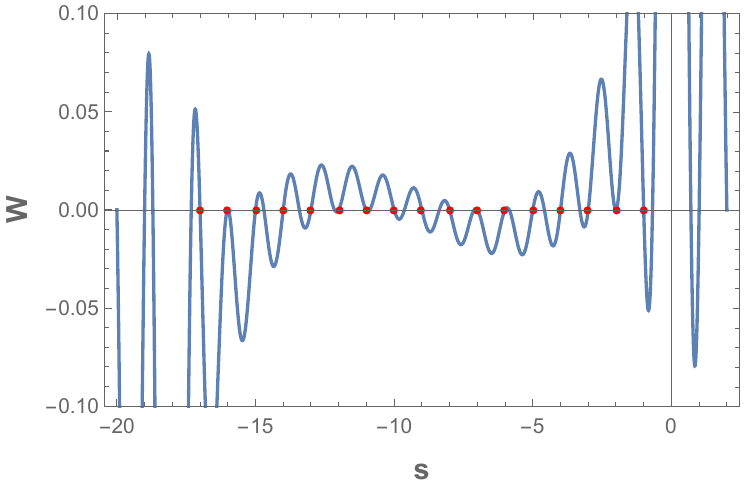}\\
(a) & (b)
\end{tabular}
\end{center}
\caption{The partition function behavior -- 
(a) the parity property  (\ref{SPF_symm}) and (b) zeros in the range $-17 \le s \le -1$
of the scalar partition $W(s,{\bf d})$ with ${\bf d}=\{2,3,6,7\}$.
The red dots in (b) correspond to the integer values of $s$.
}
\label{fig0}
\end{figure}

\begin{figure}[h!]
\begin{center}
\begin{tabular}{cc}
\includegraphics[height=4.8cm]{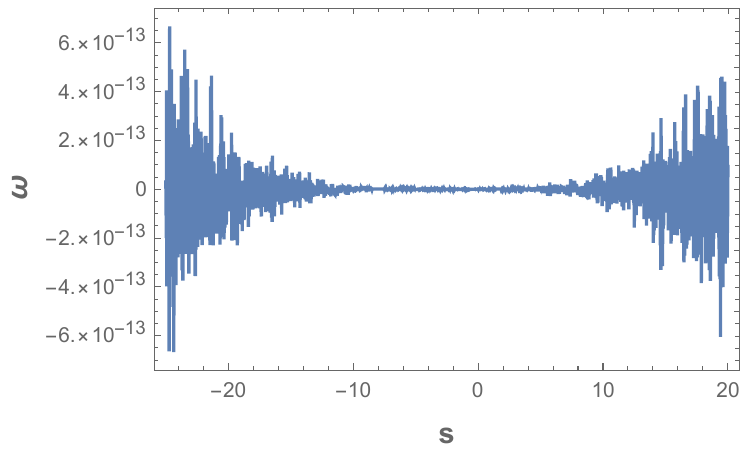}&
\includegraphics[height=4.8cm]{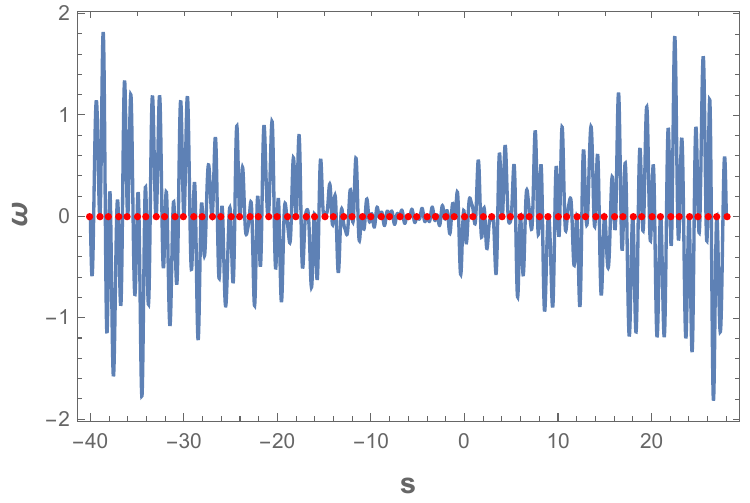}\\
(a) & (b)
\end{tabular}
\end{center}
\caption{
The behavior of the function $w(s,{\bf d},\bm \delta)$ for the unique 
generators ${\bf d}=\{2,3,6,7\}$ with  (a) $\delta_i=1$, and (b) $\bm\delta=\{1,2,1,1\}$.
The red dots in (b) correspond to the integer values of $s$.
}
\label{fig1}
\end{figure}

Consider the scalar partition $W(s,{\bf d})$ with the set of four unique elements ${\bf d}=\{2,3,6,7\}$;
the behavior of the continuous version of this partition is shown below -- 
its parity property in Fig.~\ref{fig0}(a) and zeros at negative argument values in Fig.~\ref{fig0}(b).
Derive 
a linear relation involving $W(s,{\bf d})$ 
by setting all $\delta_i$ to unity as shown in 
(\ref{matr_hatD_unique}). The elements of the vectors ${\bf d}_i$ are computed in 
(\ref{part_scalar1})  and are given by
\be
{\bf d}_1 = \{2,1,4,5\},\ 
{\bf d}_2 = \{3,-1,3,4\},\
{\bf d}_3 = \{6,-4,-3,1\},\
{\bf d}_4 = \{7,-5,-4,-1\}.
\label{unique_delta_eq_1}
\ee
The choice $\bm \delta = {\bf 1}$ is a special one as it should produce $w(s,{\bf d},\bm \delta) = 0$ everywhere
as discussed above.
The behavior of the function $w(s,{\bf d},\bm \delta)$ shown in 
Fig.~\ref{fig1}(a) confirms this prediction with high accuracy. 

The expressions for the polynomial parts of the summands in $w(s,{\bf d},\bm \delta)$
read
\bea
&&W_1(s,{\bf d}) = \frac{s^3}{1512} + \frac{s^2}{56} + \frac{437s}{3024} + \frac{113}{336}, 
\quad 
W_1(s_1,{\bf d}_1) = \frac{s^3}{240} + \frac{3s^2}{40} + \frac{193s}{480} + \frac{49}{80}, 
\nonumber  \\
&&W_1(s_2,{\bf d}_2) = -\frac{s^3}{216} - \frac{s^2}{16} - \frac{13s}{54} - \frac{23}{96}, 
\quad
W_1(s_3,{\bf d}_3) = \frac{s^3}{432} - \frac{31s}{864}, 
\nonumber \\ 
&&W_1(s_4,{\bf d}_4) =- \frac{s^3}{840} + \frac{3s^2}{560} + \frac{2s}{105} - \frac{41}{1120},
\label{unique_delta_eq_1_W1}
\eea
and the direct computation confirms the formula (\ref{poly_lin_rel}) for $w_1(s,{\bf d},\bm \delta)$.
Retaining the leading terms in (\ref{unique_delta_eq_1_W1}) we verify the result (\ref{poly_lead_terms_rel}) 
$$
(2\cdot3\cdot6\cdot7)\times
\left(\frac{1}{
2\cdot4\cdot5}-
\frac{1}{
3\cdot3\cdot4}+
\frac{1}{
3\cdot4\cdot6}-
\frac{1}{
4\cdot5\cdot7}
\right) = 
\frac{63}{10}-
\frac{7}{1}+
\frac{7}{2}-
\frac{9}{5}
= 1,
$$
and computing
$
\sigma_1({\bf d}) =18,\sigma_2({\bf d}) =98,
\sigma_1({\bf d}_i)= \{12,9,0,-3\},
\sigma_2({\bf d}_i)= \{46,35,62,91\},
$
find
$$
14\cdot\left(
\frac{3}{10}-\frac{1}{4}+\frac{3}{140}
\right) = 1,
\quad
\frac{126}{437}\cdot\left(
\frac{193}{20}-\frac{52}{9}-\frac{31}{36}+\frac{16}{35}
\right) = 1.
$$

It is instructive to consider another set of $\delta_i$ values for the same set ${\bf d}=\{2,3,6,7\}$ 
and we choose $\bm \delta=\{1,2,1,1\}$ that satisfies all conditions discussed in Section \ref{SylvCayley}.
The vectors ${\bf d}_i$ in this case read
\be
{\bf d}_1 = \{2,-1,4,5\},\ 
{\bf d}_2 = \{3,11,9,1\},\
{\bf d}_3 = \{6,-4,-9,1\},\
{\bf d}_4 = \{7,-5,-11,-1\}.
\label{unique_delta_ne_1}
\ee
The behavior of the function $w(s,{\bf d},\bm \delta)$ is shown in 
Fig.~\ref{fig1}(b) and we observe that $w(s,{\bf d},\bm \delta)=0$ at all integer points as expected
and at a countable number of real points -- which is strikingly different from the case $\bm \delta = {\bf 1}$.
The polynomial parts of the summands $W_1(s_i,{\bf d}_i)$ in $w(s,{\bf d},\bm \delta)$ are given by
\bea
&&W_1(s_1,{\bf d}_1) = -\frac{s^3}{240} - \frac{s^2}{16} - \frac{127s}{480} - \frac{9}{32}, 
\quad
W_1(s_2,{\bf d}_2) = \frac{4s^3}{891} + \frac{8s^2}{99} + \frac{379s}{891} + \frac{182}{297}, 
\label{unique_delta_ne_1_W1} \\ 
&&W_1(s_3,{\bf d}_3) = \frac{s^3}{1296} - \frac{s^2}{144} - \frac{13s}{2592} + \frac{49}{864}, 
\quad
W_1(s_4,{\bf d}_4) = - \frac{s^3}{2310} + \frac{s^2}{154} - \frac{13s}{1155} - \frac{4}{77},
\nonumber 
\eea 
validating the relation (\ref{poly_lin_rel}). 
The equalities (\ref{poly_lead_terms_rel}) turn into
$$
(2\cdot3\cdot6\cdot7)\times
\left(-\frac{1}{
2\cdot4\cdot5}+
\frac{2^3}{
3\cdot9\cdot11}+
\frac{1}{
9\cdot4\cdot6}-
\frac{1}{
11\cdot5\cdot7}
\right) = 
-\frac{63}{10}+
\frac{224}{33}+
\frac{7}{6}-
\frac{36}{55}
= 1,
$$
and
$$
14\cdot\left(
\frac{1}{4}+\frac{8\cdot 2^2}{99}-\frac{1}{36}+\frac{2}{77}
\right) = 1.
$$

\subsection{Duplicate generators}
\label{Nonunique}
\begin{figure}[h]
\begin{center}
\includegraphics[height=5.0cm]{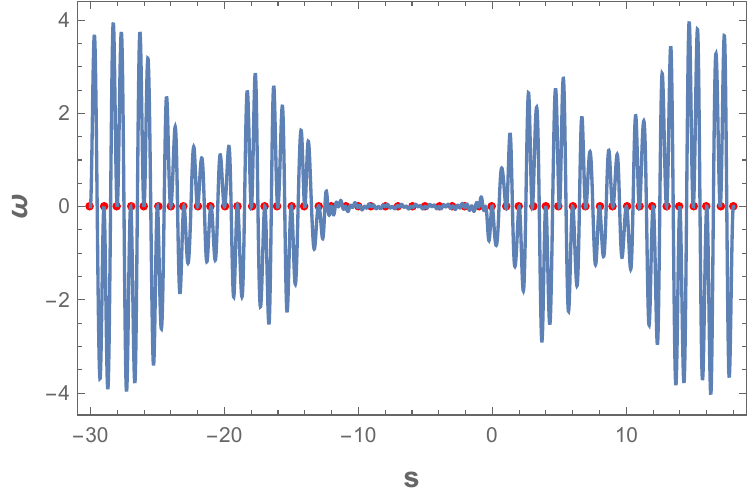}\\
\end{center}
\caption{
The behavior of the function $w(s,{\bf d},\bm \delta)$ for the set with duplicate 
generators ${\bf d}=\{2,2,5,7\}$ and $\bm\delta=\{1,3,1,1\}$.
The red dots correspond to the integer values of $s$.
}
\label{fig2}
\end{figure}

Consider an example of the scalar partitions linear relation for the sets ${\bf d}=\{2,2,5,7\}$ 
with two duplicate elements and $\bm \delta=\{1,3,1,1\}$. 
The elements of the vectors ${\bf d}_i$ are computed in 
(\ref{part_scalar1})  and equal to
\be
{\bf d}_1 = \{2,-4,3,5\},\ 
{\bf d}_2 = \{2,4,13,19\},\
{\bf d}_3 = \{5,-3,-13,2\},\
{\bf d}_4 = \{7,-5,-19,-2\}.
\label{nonunique_delta_ne_1}
\ee
The behavior of the function $w(s,{\bf d},\bm \delta)$ is shown in 
Fig.~\ref{fig2} and we observe that it oscillates around zero but indeed vanishes at integer $s$ values.
The polynomial parts of the summands in $w(s,{\bf d},\bm \delta)$ are given by
\bea
&&W_1(s,{\bf d}) = \frac{s^3}{840} + \frac{s^2}{35} + \frac{49s}{240} + \frac{29}{70}, 
\quad 
W_1(s_1,{\bf d}_1) = -\frac{s^3}{720} - \frac{s^2}{80} - \frac{3s}{160} + \frac{3}{160}, 
\nonumber \\ 
&&W_1(s_2,{\bf d}_2) = \frac{9s^3}{3952} + \frac{9s^2}{208} + \frac{1891s}{7904} + \frac{149}{416}, 
\quad 
W_1(s_3,{\bf d}_3) = \frac{s^3}{2340} - \frac{3s^2}{520} + \frac{s}{260} + \frac{63}{1040}, 
\nonumber \\ 
&&W_1(s_4,{\bf d}_4) = - \frac{s^3}{7980} + \frac{s^2}{280} - \frac{23s}{1140} - \frac{13}{560},
\label{nonunique_delta_ne_1_W1} 
\eea
verifying again the relation (\ref{poly_lin_rel}).

\section{Conclusion}
In conclusion, we present the algorithm allowing to derive infinite number of linear relations
for the scalar partitions satisfied at all {\it integer} values of the argument. 
These relations give rise to more specialized formulae (valid at the real argument values too) for
the polynomial parts of these partitions that in its turn produce a new class of relations for the 
Bernoulli polynomials of higher order. 
The general expression for leading term of the scalar partition polynomial part 
leads to a novel class of relations for the Bell polynomials.


\end{document}